\documentclass[a4paper,11pt]{article}
 \pdfoutput=1
\usepackage{color}

\usepackage{url}
\usepackage[T1]{fontenc}
\usepackage{lmodern}
\usepackage{graphicx}  
\usepackage{a4}
\usepackage{amssymb,latexsym}
\usepackage[mathscr]{eucal}
\usepackage{cite}
\usepackage[active]{srcltx}
\usepackage{fancybox}
\usepackage{ntheorem}
\usepackage{ulem}

\usepackage{cancel}

\newcommand{\referee}{\mnote{\blue{corrected, as pointed out by the referee}}}

\newcommand{\ptcheck}[1]
{\ptc{checked on #1}}

\newcommand{\syncx}[1]{\ptcxx{here starts a sync, or possibly a restrict environment}}
\renewcommand{\syncx}[1]{\ptc{syncx command here}{\color{red}#1}}

\newcommand{\ptcxx}[1]{\mnote{{\bf ptc:} {\color{red} #1}}}

\newcommand{\mnotex}[1]
{\protect{\stepcounter{mnotecount}}$^{\mbox{\footnotesize
$
\bullet$\themnotecount}}$ \marginpar{
\raggedright\tiny\em
$\!\!\!\!\!\!\,\bullet$\themnotecount: #1} }

\newcommand{\viennan}[1]{}

\newcommand{\jamesx}[1]{}
\renewcommand{\jamesx}[1]{{\mnote{{\color{blue}{\bf jg:}
#1} }}}

\newcommand{\cref}[1]{\mbox{{\color{red}FIXME; what is the 4 for?}}4\emph{\ref{#1})}}

\newcommand{\h}[2]{#1\dotfill\ #2\\\ptc{fixme}}

\newcommand{\Og}{{O}}

\newcommand{\rotpar}{\alpha}

\newcommand{\uPhi }{u_\Phihere}
\newcommand{\rPhi }{r_\Phihere}
\newcommand{\xPhi }{x_\Phihere}
\newcommand{\OmegaPhi}{\Phihere(\Omega)}

\newcommand{\hg}{{\hat g{}}}

\newcommand{\blue}[1]{{\color{blue}{#1}}}

\newcommand{\funnyr}{{r}}

\newcommand{\Lpsi}{L^2_{\psi}}

\newcommand{\Lpsione}{\zH^1_{\phi,\psi}}

\newcommand{\Lpsikg}[2]{\zH^{#1}_{\phi,\psi}(#2)}

\newcommand{\zHkpp}{\zHk_{\phi,\psi}}
\newcommand{\Hkpp}{H^k_{\phi,\psi}}
\newcommand{\zHk}{\zH^k}
\newcommand{\zH}{\mathring{H}}

\newcommand{\eqref}[1]{\eq{#1}}

\newcommand{\hs}{\cH_{\mbox{\scriptsize sing}}}

\newcommand{\beadl}[1]{\begin{deqarr}\label{#1}}
\newcommand{\eeadl}[1]{\arrlabel{#1}\end{deqarr}}%

\def \Nat{\mathbb{N}}
\newcommand{\ol}{\overline}

\def\nz{\ifmmode {I\hskip -3pt N} \else {\hbox {$I\hskip -3pt N$}}\fi}
\def\zz{\ifmmode {Z\hskip -4.8pt Z} \else
       {\hbox {$Z\hskip -4.8pt Z$}}\fi}
\def\qz{\ifmmode {Q\hskip -5.0pt\vrule height6.0pt depth 0pt
       \hskip 6pt} \else {\hbox
       {$Q\hskip -5.0pt\vrule height6.0pt depth 0pt\hskip 6pt$}}\fi}
\def\rz{\ifmmode {I\hskip -3pt R} \else {\hbox {$I\hskip -3pt R$}}\fi}
\def\cz{\ifmmode {C\hskip -4.8pt\vrule height5.8pt\hskip 6.3pt} \else
       {\hbox {$C\hskip -4.8pt\vrule height5.8pt\hskip 6.3pt$}}\fi}
\def\au{{\setbox0=\hbox{\lower1.36775ex\hbox{''}\kern-.05em}\dp0=.36775ex\hs
kip0pt\box0}}
\def\ao{{}\kern-.10em\hbox{``}}

\newcommand\Gregbeq{\begin{eqnarray}}
\newcommand\Gregeeq{\end{eqnarray}}

\def\cH{{\cal H}}

\def\h1{{\hat 1}}
\def\h2{{\hat 2}}

\def\3f{\frac{3}{2}}

\newcommand{\roscoff}[1]{}

{\catcode `\@=11 \global\let\AddToReset=\@addtoreset}
\AddToReset{figure}{section}

\DeclareFontFamily{OT1}{rsfs}{}
\DeclareFontShape{OT1}{rsfs}{m}{n}{ <-7> rsfs5 <7-10> rsfs7 <10-> rsfs10}{}
\DeclareMathAlphabet{\mycal}{OT1}{rsfs}{m}{n}

{\catcode `\@=11 \global\let\AddToReset=\@addtoreset}
\AddToReset{equation}{section}

\newcounter{mnotecount}[section]

\renewcommand{\themnotecount}{\thesection.\arabic{mnotecount}}

\newcommand{\oversetty}[2]{%
\mathop{#2}\limits^{\vbox to -.1ex{%
\kern -1.5ex\hbox{$\scriptstyle #1$}\vss}}}

\newcommand{\jlcax}[1]{}
\newcommand{\eean}{\nonumber\end{eqnarray}}

\newcommand{\kk}[1]{}

\newcommand{\beq}{\begin{equation}}

\newcommand{\FS}       
                  {F}

\newcommand{\HS} 
       {H_{\mbox{\scriptsize volume}}}

{\ptc{this should be removed in the oberwolfach version}}%

\newcommand{\eeal}[1]{\label{#1}\end{eqnarray}}
\newcommand{\bed}{\begin{deqarr}}
\newcommand{\eed}{\end{deqarr}}
\newcommand{\bedl}[1]{\begin{deqarr}\label{#1}}
\newcommand{\eedl}[2]{\arrlabel{#1}\label{#2}\end{deqarr}}

\newcommand{\loc}{\textrm{\scriptsize\upshape loc}}
\newcommand{\zmcH}{\,\,\,\,\mathring{\!\!\!\!\mycal H}}

\newcommand{\bel}[1]{\begin{equation}\label{#1}}
\newcommand{\bea}{\begin{eqnarray}}
\newcommand{\bean}{\begin{eqnarray}\nonumber}
\newcommand{\beal}[1]{\begin{eqnarray}\label{#1}}
\newcommand{\eea}{\end{eqnarray}}

\newcommand{\nn}{\nonumber}
\newcommand{\Eq}[1]{Equation~\eq{#1}}

\newcommand{\qed}{\hfill $\Box$}
\newcommand{\qedskip}{\hfill $\Box$\medskip}
\newcommand{\proof}{\noindent {\sc Proof:\ }}
\newcommand{\be}{\begin{equation}}
\newcommand{\eeq}{\end{equation}}
\newcommand{\ee}{\end{equation}}
\newcommand{\beqa}{\begin{eqnarray}}
\newcommand{\eeqa}{\end{eqnarray}}
\newcommand{\beqan}{\begin{eqnarray*}}
\newcommand{\eeqan}{\end{eqnarray*}}
\newcommand{\ba}{\begin{array}}
\newcommand{\ea}{\end{array}}

\newcommand{\hyp}{\mycal S}

\newcommand{\mnote}[1]
{\protect{\stepcounter{mnotecount}}$^{\mbox{\footnotesize
$
\bullet$\themnotecount}}$ \marginpar{
\raggedright\tiny\em
$\!\!\!\!\!\!\,\bullet$\themnotecount: #1} }

\newcommand{\warn}[1]
{\protect{\stepcounter{mnotecount}}$^{\mbox{\footnotesize
$
\bullet$\themnotecount}}$ \marginpar{
\raggedright\tiny\em
$\!\!\!\!\!\!\,\bullet$\themnotecount: {\bf Warning:} #1} }

\newcommand{\R}{\mathbb R}

\newcommand{\N}{\mathbb N}

\newcommand{\eq}[1]{(\ref{#1})}

\newcommand{\ptc}[1]{\mnote{{\bf ptc:}#1}}

\newcommand{\Ric}{\mbox{\rm Ric}}

\newcommand{\beqar}{\begin{deqarr}}
\newcommand{\eeqar}{\end{deqarr}}

\newcommand{\beaa}{\begin{eqnarray*}}
\newcommand{\eeaa}{\end{eqnarray*}}

\newcommand{\bethm}{\begin{theorem}}
\newcommand{\et}{\end{theorem}}
\newcommand{\bl}{\begin{Lemma}}

\newtheorem{Theorem} {\sc  Theorem\rm} [section]
\newtheorem{theorem} [Theorem] {\sc  Theorem\rm}
\newtheorem{theo} [Theorem] {\sc  Theorem\rm}

\newtheorem{Corollary} [Theorem] {\sc  Corollary\rm}

\newtheorem{Lemma} [Theorem] {\sc  Lemma\rm}
\newtheorem{Proposition} [Theorem] {\sc  Proposition\rm}

\theorembodyfont{\upshape}

\newtheorem{remark}[Theorem]{\sc Remark\rm}

\theoremsymbol{\ensuremath{\diamondsuit}}
\newcommand{\fcoco}{\small}
\theorembodyfont{\fcoco}\theoremseparator{}\theoremindent0.5cm

\theoremstyle{nonumberplain}\theorembodyfont{\fcoco}
\theoremseparator{}\theoremindent0.5cm

\DeclareFontFamily{OT1}{rsfs}{}
\DeclareFontShape{OT1}{rsfs}{m}{n}{ <-7> rsfs5 <7-10> rsfs7 <10-> rsfs10}{}
\DeclareMathAlphabet{\mycal}{OT1}{rsfs}{m}{n}

{\catcode `\@=11 \global\let\AddToReset=\@addtoreset}
\AddToReset{figure}{section}

{\catcode `\@=11 \global\let\AddToReset=\@addtoreset}
\AddToReset{equation}{section}

\renewcommand{\Og}{O}
\renewcommand{\blue}[1]{#1}
\renewcommand{\referee}{}

\begin{document}

\title{On Carlotto-Schoen-type scalar-curvature gluings\footnote{Preprint UWThPh 2016-19}}

\author
{Piotr T.~Chru\'sciel,
%
%
%
%
Erwann Delay}
%
%
%

\maketitle
\begin{abstract}
 We  carry out a Carlotto-Schoen-type  gluing  with interpolating scalar curvature on cone-like sets, or deformations thereof, in the category of smooth Riemannian asymptotically Euclidean metrics.
\end{abstract}

\section{Introduction}
 \label{s23III13.1}

In an outstanding paper~\cite{CarlottoSchoen} Carlotto and Schoen have constructed non-trivial asymptotically Euclidean scalar-flat metrics which are Minkowskian outside of a solid cone. The reading of their paper suggests strongly, though no explicit statements are made, that their construction can be generalised as follows:

\begin{enumerate}
 \item Rather than gluing an asymptotically Euclidean metric to a flat one, any two asymptotically Euclidean metrics $g_1$ and $g_2$ can be glued together.
 \item In the spirit of~\cite{ErwannInterpolating}, the gluings at zero-scalar curvature can be replaced by gluings where the scalar curvature of the final metric equals
     $$
      \chi R(g_1) + (1-\chi) R(g_2)
     $$
     with a cut-off function $\chi$.
      \item The geometry of the gluing region can be somewhat more general than the interface between two rotationally-symmetric cones.
\end{enumerate}

The aim of this paper is to give detailed proofs of the above. While the overall structure of the argument is rather similar, the details are different:
Both the proof in \cite{CarlottoSchoen} and ours require a weighted Poincar\'e inequality, with nonstandard weights adapted to the geometry of the problem at hand. We provide a version of this inequality different from the one in~\cite{CarlottoSchoen}, which might be of more general interest. We also provide a detailed proof of the gluing with exponential weights near the boundary, only sketched in~\cite{CarlottoSchoen}.
Finally, we replace some of the tailor-made arguments in \cite{CarlottoSchoen} by off-the-shelf results  from \cite{ChDelay}.

\begin{figure}
\vspace{-.9cm}
  \centering
    \includegraphics[scale=.2]{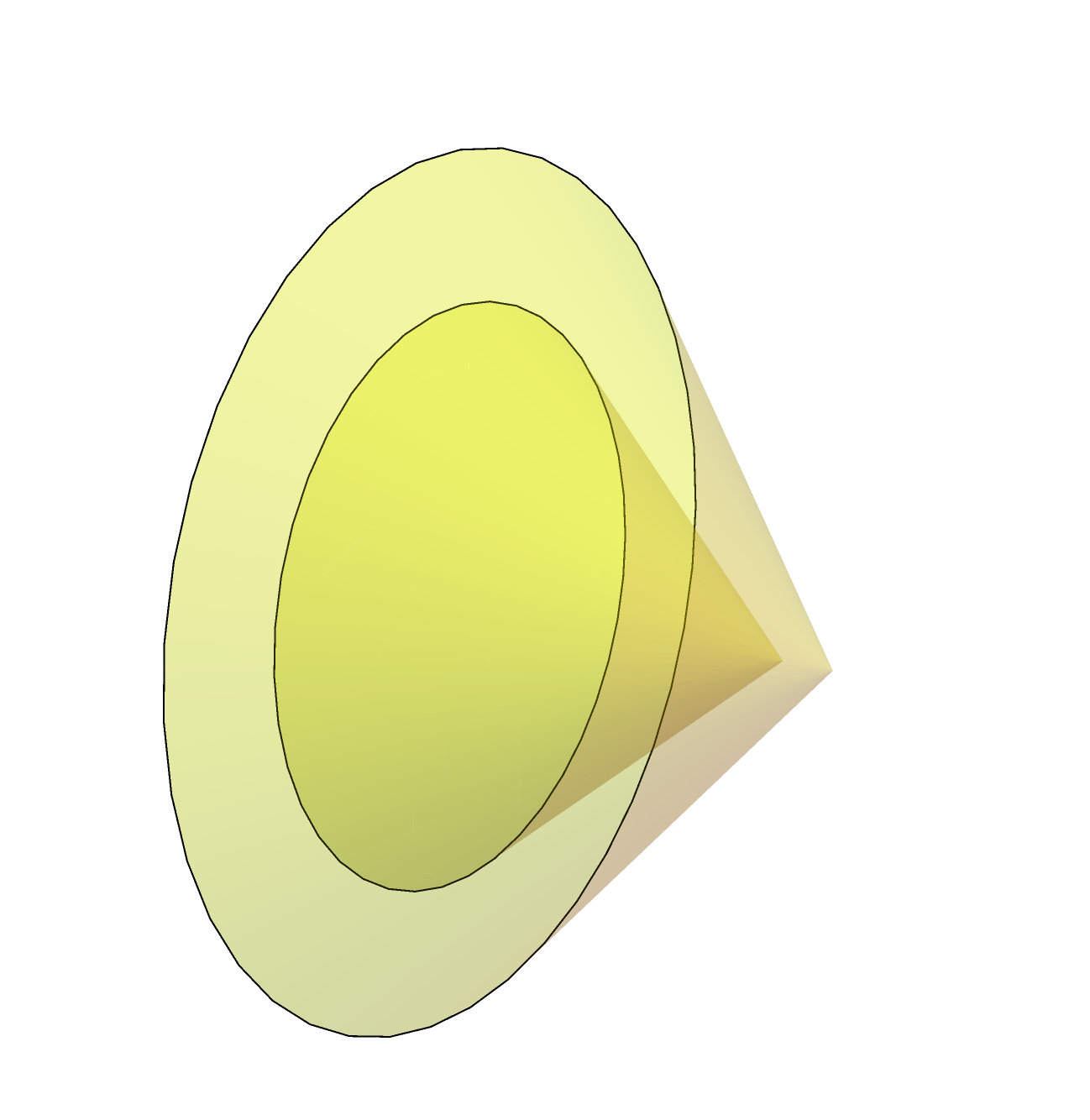}
    \includegraphics[scale=.25]{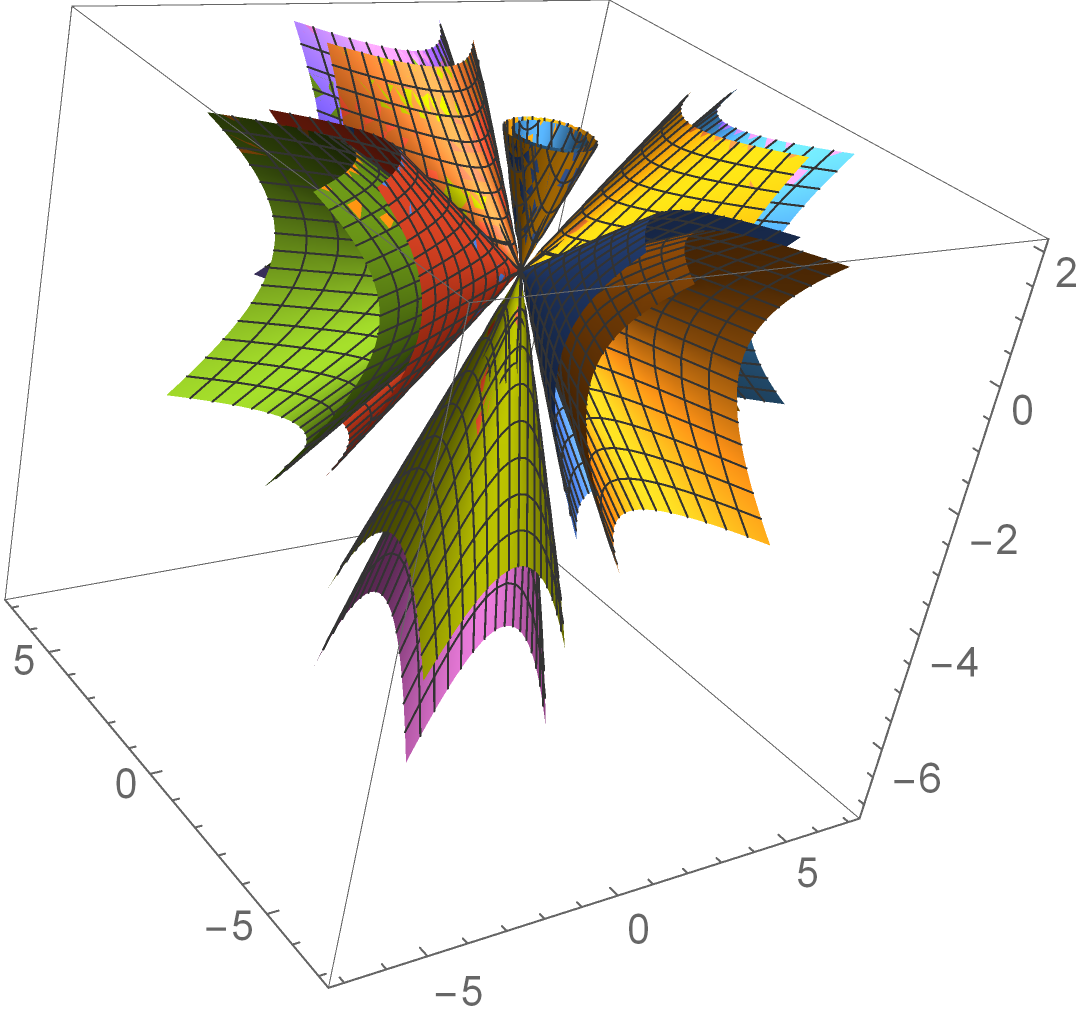}
    \vspace{-.4cm}
\caption{Left picture: If the second metric is flat and both metrics are scalar-flat, then the new metric is flat  inside the smaller cone,  scalar-flat everywhere, and coincides with the first one  outside the larger cone. Both cones extend to infinity, and their tips are located very far in the asymptotically Euclidean  region.
Right picture: Iterating the construction, one can embed any finite number of distinct initial data sets into Minkowskian data, or paste-in Minkowskian data inside several cones into a given data set.
    \label{F11IX16.1}}
\end{figure}

A special case of our main Theorem~\ref{ThefulltheoremAE} below provides the following variant of the Riemannian-geometry version of the main theorem of~\cite{CarlottoSchoen}:

\begin{theo}
  \label{T11IX16.1}
Let $n\ge 3$. Consider two smooth $n$-dimensional asymptotically Euclidean Riemannian metrics $g_1$ and $g_2$ and two nested cones with non-zero opening angles and vertices displaced along the symmetry axis, see Figure~\ref{F11IX16.1}. Simultaneously scaling-up the cones if necessary, or simultaneously shifting them to the asymptotic region, there exists a smooth asymptotically Euclidean  metric $g$ which coincides with $g_1$ outside of the larger cone and with $g_2$ inside the smaller one, with the scalar curvature $R(g)$ lying between $R(g_1)$ and $R(g_2)$ in the intermediate region.
\end{theo}

In particular if $R(g_1)=R(g_2)=0$, then $R(g)=0$.

In Section~\ref{S7X16.1} we show how Theorem~\ref{T11IX16.1} follows from our main Theorem~\ref{ThefulltheoremAE} below.
The key to the proof of this last theorem is a weighted Poincar\'e inequality involving radially-scaled exponential weights, proved in Proposition~\ref{PL4VIII14.2} below.
The rest of the proof is a verification of the hypotheses of the rather general results proved in~\cite{ChDelay}.

As already mentioned,  Theorem~\ref{ThefulltheoremAE} below allows more general cones than the ones in~\cite{CarlottoSchoen}, as described at the beginning of Section~\ref{S4X16.1}. In fact, our arguments apply to a large class of deformed cones as well, e.g.\ ``logarithmically-rotated ones'', cf.\ Theorem~\ref{T18X16.1} and the discussion in Section~\ref{s19X16.1}.

\section{Definitions, notations and conventions}
 \label{sec:def}

We use the summation convention throughout, indices are lowered with $g_{ij}$
and raised with  its inverse $g^{ij}$.

We will have to frequently control the asymptotic behavior of the objects at hand. Given a tensor field $T$ and a function $f$,  we will write
$$
T=\Og (f),
$$
when there exists a constant $C$ such that the $g$-norm of $T$ is dominated by $Cf$.

A metric $g$ on a manifold $M$ will be said to be asymptotically Euclidean (AE) if $M$ contains a region, diffeomorphic to the complement of a ball in $\R^n$, on which the metric $g$ approaches the Euclidean metric $\delta$ as one recedes to infinity.

Let $\phi$ and $\psi$ be two smooth strictly
positive functions on an $n$-dimensional manifold
$M$. The function $\psi$ will be used to control the growth of the fields involved near boundaries or in the asymptotic regions, while $\phi$ will control how the growth is affected by derivations.
For $k\in \Nat$ let $\Hkpp(g) $ be the space of $H^k_\loc$
functions or tensor fields such that the norm\footnote{The reader
is referred to~\cite{Aubin,Aubin76,Hebey} for a discussion of
Sobolev spaces on Riemannian manifolds.}
\be \label{defHn}
 \|u\|_{\Hkpp (g)}:=
 (\int_M(\sum_{i=0}^k \phi^{2i}|\nabla^{(i)}
 u|^2_g)\psi^2 d\mu_g )^{\frac{1}{2}}
\ee
is finite, where
$\nabla^{(i)}$ stands for the tensor $\underbrace{\nabla ...\nabla
}_{i \mbox{ \scriptsize times}}u$, with $\nabla$ being the
Levi-Civita covariant derivative of $g$. We assume throughout that
the metric is at least $W^{1,\infty}_\loc$; higher
differentiability will be usually indicated whenever needed.

For
$k\in \Nat$ we denote by $\zHkpp $ the closure in $\Hkpp$ of the
space of smooth functions or tensors which are compactly
 supported in $M$, with the norm induced from
$\Hkpp$.
The $\zHkpp $'s are Hilbert spaces with the obvious scalar product
associated with the norm \eq{defHn}. We will also use the following
notation
$$
\quad \zHk  :=\zHk  _{1,1}\,,\quad
L^2_{\psi}:=\zH^0_{1,\psi}=H^0_{1,\psi}\,,
$$
so that $L^2\equiv \zH^0:=\zH^0_{1,1}$.

\blue{
\referee
For $\phi$ and $\varphi$  --- smooth strictly positive functions
on $M$, and for $k\in\N$ and $\alpha\in [0,1]$, we define
$C^{k,\alpha}_{\phi,\varphi}$ to be the space of $C^{k,\alpha}$
functions or tensor fields  for which the norm
$$
\begin{array}{l}
\|u\|_{C^{k,\alpha}_{\phi,\varphi}(g)}=\sup_{x\in
M}\sum_{i=0}^k\Big(
\|\varphi \phi^i \nabla^{(i)}u(x)\|_g\\
 \hspace{3cm}+\sup_{0\ne d_g(x,y)\le \phi(x)/2}\varphi(x) \phi^{i+\alpha}(x)\frac{\|
\nabla^{(i)}u(x)-\nabla^{(i)}u(y)\|_g}{d^\alpha_g(x,y)}\Big)
\end{array}$$ is finite. Here $d_g(p,q)$ denotes the $g$-distance between $p$ and $q$. Note that the last term is not needed when $\alpha=0$.
}

Let $g$, $g_0$ be two Riemannian metrics and let $W$ be a space of two-covariant symmetric tensors. We will write $g\in M_{g_0+W} $ if $g-g_0 \in W$. An important example is provided by metrics
$g\in M_{\delta+C^{k+4}_{r,r^{\epsilon}}}$, where $\delta$ denotes a metric which is Euclidean in the asymptotic end, and $r$ is a (strictly) positive smooth function which equals $|\vec x|$ in the explicitly Euclidean coordinates in the asymptotic region.
 Such metrics are thus  \emph{asymptotically Euclidean}, with the difference $g-\delta$ decaying as $O(r^{-\epsilon})$, and with derivatives of order $1\le\ell\le k+4$ decaying  as $O(r^{-\epsilon-\ell})$.

\section{Gluing on sets which are scale-invariant at large distances}
 \label{S4X16.1}

Let $S(p,R)\subset \R^n$ denote a sphere of radius $R$ centred at $p$.

Let $\Omega \subset \R^n$ be a domain with smooth boundary, thus $\Omega$ is open and connected, and $\partial \Omega =\overline \Omega \setminus \Omega$ is a smooth manifold. In other words, $
\overline \Omega$ is a smooth (connected) submanifold of $\R^n$ with smooth boundary and $\Omega$ is its interior.
We further assume that $\partial \Omega$ has \emph{exactly two} connected components, and that
there exists $R_0\ge 1$ such that
\blue{\referee
 the boundary $\partial \Omega_S$ of}
\bel{7X16.24}
 \Omega_{S}:= \Omega \cap S(0,R_0)
\ee
also has again exactly two
 connected components, with
\bel{7X16.25}
 \Omega \setminus B(0,R_0)=\{\lambda p \ | \ p \in \Omega_{S}\,, \lambda \ge 1\}
 \,.
\ee
When \eq{7X16.25} holds, we say that $\Omega$ is \emph{invariant under scaling at large distances}.

We let $x:\overline \Omega\to \R$ be any smooth defining function for $\partial \Omega$ which has been chosen so that
\bel{14VI14.1}
 \mbox{$x(\lambda p) = \lambda x(p)$ for $\lambda \ge 1$ and for  $p, \lambda p \in \Omega\setminus B(0,R_0)$.}
\ee
Equivalently, for $p\in\Omega_S$ and $\lambda$ larger than one, we require $x(\lambda p) = \lambda x_S(p)$, where $x_S$ is a defining function for  $\partial\Omega_S$ within $S(0,R_0)$.

By definition of $\Omega$ there exists a constant $c$ such that the distance function $d(p)$ from a point $p\in\Omega$ to $\partial \Omega$ is smooth for all $d(p)\le c r(p)$
for some (perhaps small) constant $c$ and where, as before, $r=|\vec x|$ in the explicitly Euclidean coordinates in the asymptotic end.
 The function $x$ can be chosen to be equal to $d$ in that region.

A simple example, considered in~\cite{CarlottoSchoen}, satisfying the above is the following: let $0<\theta_1<\theta_2<\pi$ and let $\Omega_{\theta_1,\theta_2}$ any domain with smooth boundary such that
\bel{7X16.31}
 \Omega_{\theta_1,\theta_2} \setminus B(0,R_0)=( \mathring C_{\theta_2}\setminus C_{\theta_1}) \setminus B(0,R_0)
 \,,
\ee
where $C_\theta$ is a closed solid  cone with aperture $\theta>0$;
compare the left Figure~\ref{F11IX16.1}. Here, and elsewhere, $\mathring U$ denotes the interior of a set $U$. In this example $\Omega_S:=\Omega_{\theta_1,\theta_2}
\cap S(0,R_0)$ is connected, with a boundary which has two connected components.

\begin{figure}[th]
  \centering
    \includegraphics[scale=.3]{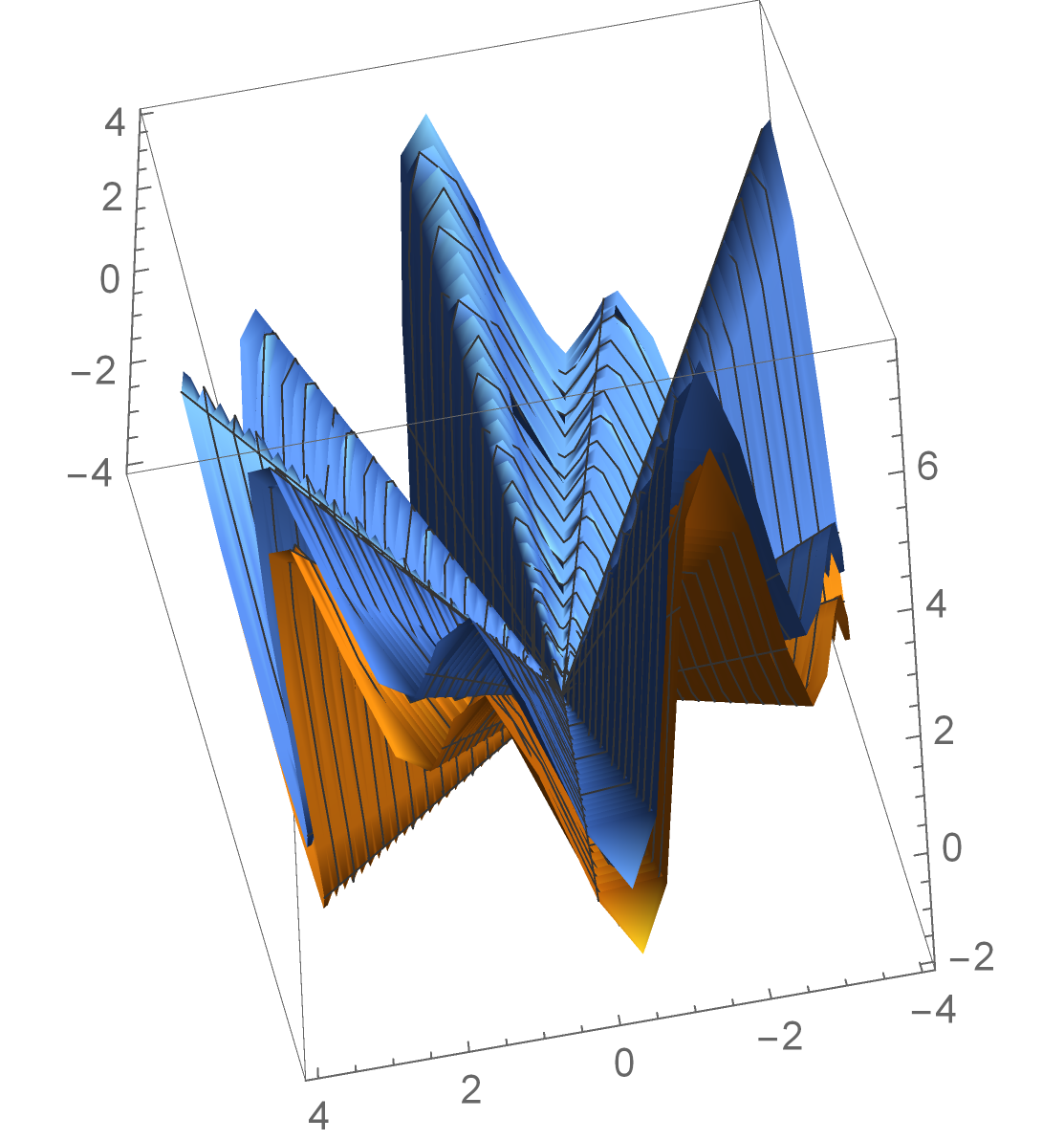}
    \includegraphics[scale=.35]{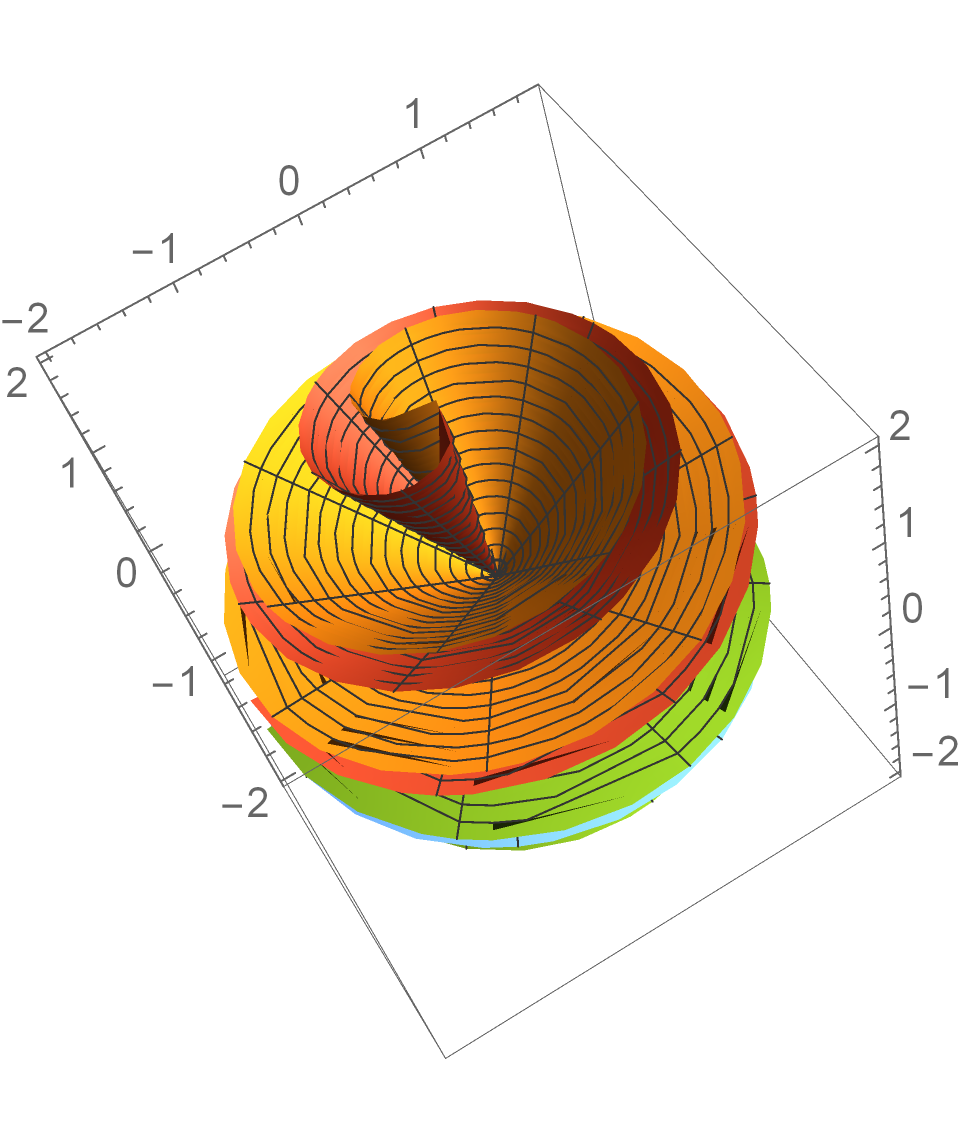}
\caption{Left figure: A possible set $\Omega$,  located between the two surfaces. The metric coincides with $g_1$ above the first surface, with $g_2$ below the second one, and is scalar-flat everywhere if both $g_1$ and $g_2$ were.
Right figure: A scale-invariant thickening of the displayed hypersurface provides a set $\Omega$ so that $\R^3\setminus \Omega$ has two components. The metric will coincide with $g_1$ in one component, and with $g_2$ in the other.
    \label{F17X16.1}}
\end{figure}

Some other examples can be found in Figure~\ref{F17X16.1}
\begin{figure}[th]
  \centering
    \includegraphics[height=4cm]{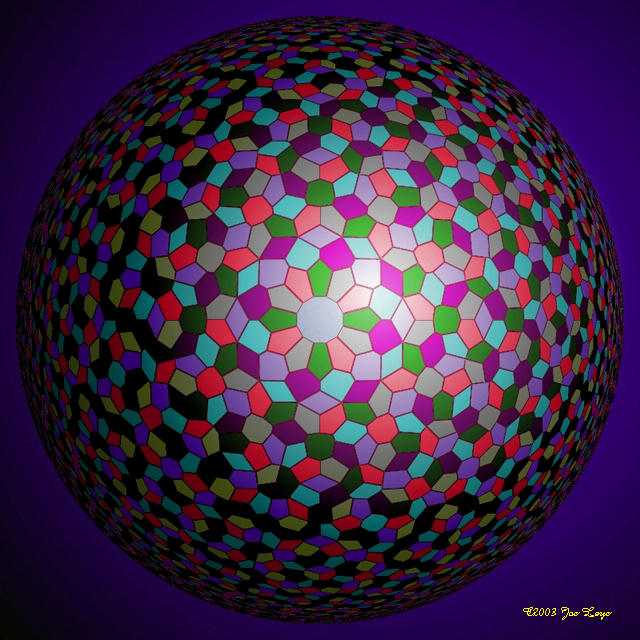}
    \includegraphics[height=4cm]{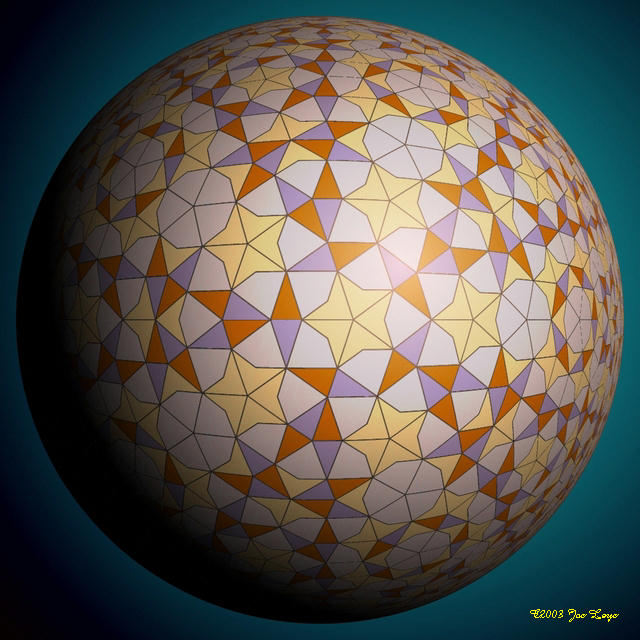}
    \includegraphics[height=4cm]{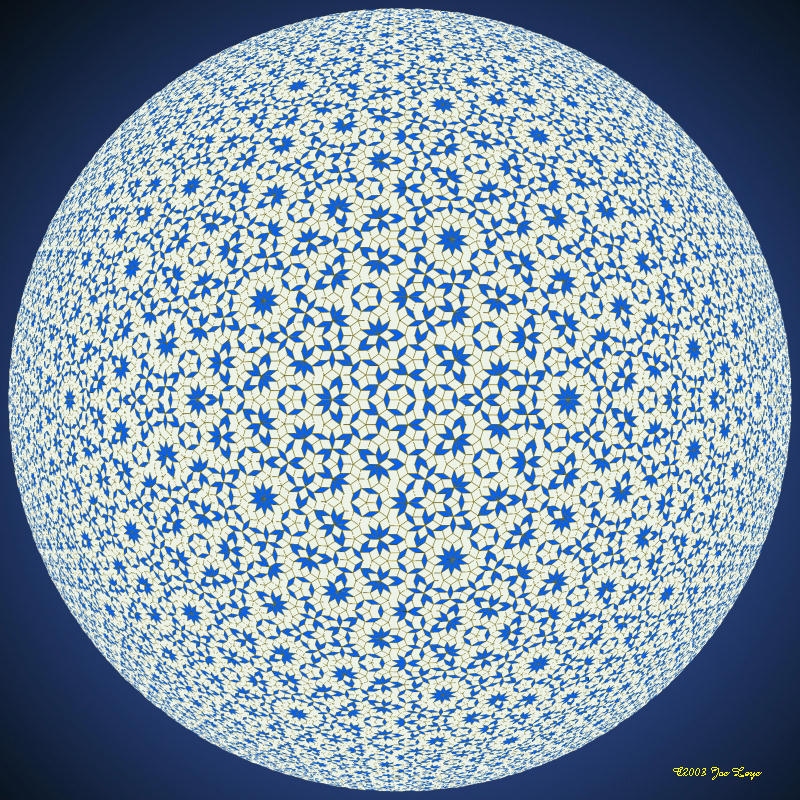}
\caption{After iterating the gluing a finite number of times (first two pictures), or an infinite number of times (last picture) one can obtain aesthetically pleasing sets $\Omega_S$. The infinite iteration might require  going to larger and larger distance for each individual gluing.  Images by Jos Leys, with kind permission of the author.
    \label{F17X16.1+}}
\end{figure}

We will denote by $r$  a  smooth positive function which coincides with $|\vec x|$ for $|\vec x |\ge 1$.

For $\beta, s, \mu \in \R$ we define
\bel{5X16.11}
 \phi=\left(\frac{x}{r}\right)^2r=\frac{x^2}r
 \,,
  \quad
  \psi = \funnyr^{-n/2-\beta}\left(\frac{x}{r}\right)^\sigma e^{-sr/x}=:\funnyr^{\mu}{x}^\sigma e^{-sr/x}
\ee
on $\Omega$. The factor $r^{-\beta}$ is directly related to the large-distance radial behavior of the solutions, while the exponential factor with $s>0$ will force the solutions {we construct} to vanish at all orders at $\partial \Omega$ (compare \eq{6X16.11} below).

We will make a gluing-by-interpolation of scalar curvature. The main interest is that of scalar-flat metrics, which then remain scalar-flat, or for metrics with positive scalar curvature, which then remains positive.  Since the current problem is related to the construction of initial data sets for Einstein equations, in some general-relativistic matter models, such as   Vlasov  or dust, an  interpolation of scalar curvature is  of direct interest.

 In order to carry out the interpolation, recall that $\partial \Omega$ has exactly two connected components. We denote by $\chi$ a smooth function with the following properties:
\begin{enumerate}
 \item $0\le \chi \le 1$;
   \item $\chi$ equals one in a neighborhood of one of the components and equals zero in a neighborhood of the other component;
       \item on $\Omega\setminus B(0,R_0)$ the function $\chi$ is required \blue{\referee to depend only upon the variable $p$ in \eqref{7X16.25}.}
\end{enumerate}

Starting  with two AE metrics $ g $ and $ \hat g $  we  define
\beal{1XI16.1}
 &
 g_\chi:=\chi\hat g+(1-\chi)g
 \,,
 &
\\
 &
   R_\chi:=\chi R(\hat g)+(1-\chi)
    R(g)\,,
    \quad
 \delta R_\chi:=R_\chi
 -R \big(g_\chi\big)
 \,.
\eeal{1XI16.2}

We denote by $P_{g}$ the linearisation of the scalar-curvature  operator at a metric $g$, and will write $P$ for $P_g$ when $g$ is obvious in the context. Its formal adjoint $P_g^*$ reads
\bel{7X16.1}
 P_g^* (f) = \nabla \nabla f - \Delta f g -f \Ric (g)
 \,,
\ee
where $\Delta= \nabla^i \nabla_i$.

Letting $\Omega$, $x$ and $r$ be as just described, we have:

\begin{theorem}
 \label{ThefulltheoremAE}
 Let $\epsilon>0$, $k>n/2$, $\beta\in[-(n-2),0)$,    {$\tilde \beta<-\epsilon<\beta$, $s >0$, and $\sigma \neq -1/2 $.}
%
%
{
Set
\bel{14III19.1}
 \phi=\frac{x^2}r
 \,,
  \quad
  \psi = \funnyr^{-n/2-\beta}\left(\frac{x}{r}\right)^\sigma e^{-sr/x}
  \,.
\ee
}
Suppose that $g\in M_{\delta+C^{k+4}_{r,r^{\epsilon}}}\cap C^\infty$.
For all real numbers $\sigma$ and $s>0$  and
$$
 \mbox{all smooth metrics
    $\hg $ close enough to $ g$ in $ C^{k+4}_{r,r^{- {\tilde \beta}   }}(\Omega )$}
$$
there exists on $\Omega$ a unique smooth two-covariant symmetric tensor field $ h  $ of the form
\bel{6X16.12}
      h
   = \phi^4\psi^2  P^*_{ g_\chi}( \delta N)\in \psi^2  \phi^2\Lpsikg{k+2}{g}
\ee
such that the metric $ g_\chi+ h   $ solves
\bel{fullcolle}
 R
    \left[g_\chi+ h  \right]= \chi R(\hat g) +(1-\chi)
    R(g)
 \,.
\ee
The   tensor field $  h   $ vanishes at $\partial \Omega$ and can be $C^{\infty}$-extended by zero across $\partial \Omega$, leading to a smooth  asymptotically Euclidean metric $g_\chi+ h  $.
\end{theorem}

\begin{remark}
  \label{R7X16.1}
{\rm
Some comments about the asymptotic behaviour of the metrics involved are in order. The requirement $g\in M_{\delta+C^{k+4}_{r,r^{\epsilon}}}\cap C^\infty$ (equivalently, $g-\delta \in {C^{k+4}_{r,r^{\epsilon}}}\cap C^\infty$) guarantees that $g$ is asymptotically Euclidean, with $g$ approaching the Euclidean metric as $O(r^{-\epsilon})$. Note that $\epsilon$ can be arbitrary small. Next, $g-
\hat g$ is required to fall off as $r^ {\tilde \beta}  $, so that $g_\chi$ will approach $\delta$ as
{$O(r^{\max (-\epsilon,  {\tilde \beta}  )})$}.
 The final metric $g_\chi+ h$ will a priori approach $\delta$ somewhat slower.
Indeed, the proof below shows that there exists a constant $C  $ such that
\bea
 \| h  \|_{\psi^2\phi^2\mathring H^{k+2}_{\phi,\psi}(g_\chi) }
  \leq C \left\|\delta R_\chi
        \right\|_{ \psi^2 \mathring H^{k }_{\phi,\psi}(g_\chi) }
 \,.
\eeal{estimatesolutionfullC}
{(The right-hand side is finite and small only if the weight $\beta$ is larger than $\max( {\tilde \beta}  ,-\epsilon)$, which explains the hypothesis.)}
{(Also note that $\delta N \in \mathring H^{k+4}_{\phi,\psi}(g_\chi)$, in particular $\delta N=o(r^{\beta})$ .)}

 \Eq{estimatesolutionfullC} implies (compare \eq{6X16/22} below)
\bel{6X16.11}
 h=o (r^{-n-\beta +2} (x/r)^{ \sigma + 4 - n} e^{-sr/x})
  =o (r^{-n-\beta +2})
 \,.
\ee
Since $\beta<0$, \eq{6X16.11} shows that the decay of the metric to the flat-one is \emph{most likely} slower than the Schwarzschildian decay rate $O(r^{-(n-2)})$.
}
\end{remark}

\begin{remark}
  \label{R6X16.1}
{\rm
A well defined ADM energy-momentum  requires the metric to approach the flat-one as $  o(r^{-(n-2)/2})$. Comparing with \eq{6X16.11}, this will be satisfied if in Theorem~\ref{ThefulltheoremAE} we require {
 $\epsilon > (n-2)/2$ and}
\bel{13VI14.7}
 -\frac{ n-2}2
  \le \beta \le 0
 \,.
\ee
%
%
}
\end{remark}

\begin{remark}
  \label{R6X13.3}
{\rm
An analogous result holds if $g$ is not necessarily smooth: If $g \in   M_{\delta+C^{k+4}_{r,r^{\epsilon}}} $, then the construction can be modified as in~\cite{ChDelayHilbert} so that the final metric $g_\chi+h$ will be of $C^{k+4}$ differentiability class.
}
\end{remark}

\proof
The idea of the proof is to use some general results proved in \cite{ChDelay}.
Now, the whole setup of \cite{ChDelay} was geared towards the general relativistic constraint equations, which involve both a metric $g$ and an extrinsic curvature tensor $K$. However, the analysis there applies\emph{ mutatis mutandis} to the Riemannian problem, where only the scalar curvature $R$ of the metric $g$ is considered, by setting $K\equiv 0 $, and e.g.\ setting to zero the vector field $Y$ appearing  in~\cite{ChDelay}. It then suffices to ignore all the conditions imposed in~\cite{ChDelay} on the part of the equations there which involve a vector field $Y$. For example, in the current setting the inequality (3.4) of \cite{ChDelay} should be replaced by %
\be
\label{it1kerbis} C\| {\phi^2} P_{g_0}^* (N)\|_{\Lpsi(g_0)}
 \geq
 \|N\|_{\Lpsione(g_0)}\,.
\ee
with $
 P_g^*  $ defined in \eqref{7X16.1}.

Given this, we
 start by showing that Theorem~3.6  of~\cite{ChDelay}  with $K\equiv 0 \equiv Y$ applies, with $g_0$ taken to be 
 a Ricci flat metric which is uniformly equivalent to the Euclidean one. 
This requires, first, checking the ``scaling property'' of the functions $\phi$ and $\psi$ {as defined in~\cite[Appendix~B]{ChDelay}.
}
This is a routine exercise, using scaling in $r$ and the arguments in~\cite[Appendix~B]{ChDelay}. Next, note \blue{first that with a Ricci flat metric $g_0$,} Equation~\eq{7X16.1} reads
\bel{7X16.2}
 P_{g_0}^* (f) = \nabla \nabla f - \Delta f g_0
 \,.
\ee
Since $\beta \ne 0$, the
\blue{
inequality {\eqref{it1kerbis}}
\referee
needed for \cite[Theorem~3.6]{ChDelay}
}
 follows directly from Proposition~\ref{PL4VIII14.2} below used twice, first for {$N$ and then for $\nabla N$}. See also Remark~\ref{R1XI16.1}.
 When $g_0$ is \emph{not} Ricci flat, the Ricci terms in $P^*$ are dominated by the remaining ones by choosing {$x/r$} small enough in the relevant estimates: indeed, the term
$$
\phi^2  \Ric(g) N=(x^2/r)^2 O(r^{-\epsilon -2})N=O(r^{-\epsilon}) (x/r)^4 N =o(1) N
$$
can be absorbed in the right-hand side of \eqref{it1kerbis} after using the triangle inequality on the left-hand side, resulting in \eqref{it1kerbis} with $C$ replaced by, say, $2C$.
Further note that the integrals involving $u$ in the right-hand side
of \eq{9IX14.110} will vanish if the support of $u$ lies in the region $|\vec x|\ge R$.

Now, the  operator $\pi_{{\mathcal K}_0^\perp}$ appearing in the statement of \cite[Theorem~3.6]{ChDelay} projects on the space ${\mathcal K}_0^\perp$ of  \emph{static KIDs}, namely the kernel of $P^*$. We claim that under our hypotheses this space is trivial, and so the
 operator
${L}_{\phi,\psi}$ of \cite[Theorem~3.6]{ChDelay} is an isomorphism:
Indeed, it is well-known that static KIDs that are $o(1)$ along some cone $C_\varepsilon$ have to vanish everywhere. (This is essentially~\cite[Proposition~2.1]{ChMaerten}, together with unique continuation for solutions of the KID equation.)
  So it suffices to show that functions $\delta N\in \Lpsikg{k+4}{g}$ are $o(1)$ in some cone contained in $\Omega$. For this let us choose a cone $C_\varepsilon \subset \Omega$ on which $x/r$ is bounded away from zero. Then the restriction $\delta N|_{C_\varepsilon}$ is in the space
\be\label{cHdef}
 \zmcH _k^\beta (C_\varepsilon)  := \zHk _{r,\funnyr^{-n/2-\beta}} (C_\varepsilon)
 \,.
\ee
We note the inclusions~\cite{Bartnik86}
\bel{sobinc} C^{\beta'}_k(C_\varepsilon)\subset \zmcH _k^\beta(C_\varepsilon)\,, \ \beta'<
\beta\,, \ \mbox{ and } \  \quad\zmcH _k^\beta(C_\varepsilon) \subset C^\beta_{\lfloor k-n/2
\rfloor}(C_\varepsilon)\,, \ k>n/2
 \,.
\ee
In fact~\cite{Bartnik86}
\bel{sobinc2} f\in \zmcH _k^\beta(C_\varepsilon) \,, \ k>n/2 \quad
\Longrightarrow \quad f=o(r^\beta)
 \,.
\ee
So, the requirement
that there are no KIDs in the space under consideration will be satisfied for
\bel{13VI14.1}
 \beta\le0
 \,.
\ee
(A more careful inspection shows (compare \eq{6X16/22} below) that
\bel{sobincAEexp} \delta N \in \mathring H^k_{\phi,\psi}  \,, \ k>n/2 \quad
\Longrightarrow \quad \delta N =o\Big(r^\beta\big(\frac{x}{r}\big)^{-\sigma-n}e^{sr/x}\Big)
 \,,
\ee
but a possible blow-up of $\delta N$  when $x$ tends to zero is irrelevant as long as $\delta N$ is forced to go to zero when receding to infinity away from $\partial \Omega$.)

Summarising,
\cite[Theorem~3.6]{ChDelay} applies and shows that the map $L_{\phi,\psi}$ there is an isomorphism.

We now want to use the
inverse function-type theorem
as in
 \blue{ \cite[Theorem~3.9]{ChDelay}}
\referee
%
  to solve the equation \eq{fullcolle}; equivalently
\bel{7X16.4}
 R
    \left[g_\chi+ h  \right] - R(g_\chi) = \delta R_\chi
 \,.
\ee
This requires checking differentiability of the map defined there. This will follow by standard arguments if the metrics of the form $g_\chi+h$, with $h$ given by  \eq{6X16.12}, are asymptotically Euclidean.  Now, given a perturbation $\delta R$ of the Ricci scalar on $\Omega$, the linearized perturbed
metric is obtained from the
solution $\delta N$ of the equation
\bel{23VI14.3}
 \psi^2 L_{\phi,\psi}\delta N \equiv P\underbrace{\phi^4\psi^2  P^*\delta N}_{=: h  }= \delta R\in \psi^2
\Lpsikg{k}{g}
\,.
\ee

In cones $C_\varepsilon$ staying away from the boundary we have $\delta N \in  \zmcH _{k+4}^\beta(C_\varepsilon)$, and since on $C_\varepsilon$ we have  $\psi\sim r^{-n/2-\beta}$  and $\phi \sim r$
we obtain
\bea
   h
   & = & \phi^4\psi^2  P^* \delta N
   =
    r^{-n-2\beta+4} o(r^{\beta-2})  =
   o(r^{-\beta-n+2})
  \,.
\eeal{13VI14.5}
We conclude that $h$  will decay to zero in such cones provided that
\bel{13VI14.6}
  -(n-2) \le \beta
 \,.
\ee
A more careful treatment, without invoking interior cones, uses the weighted Sobolev embedding
\blue{
\referee
\bean
 h\in \psi^2\phi^2\mathring H^{k+2}_{\phi,\psi}(g_\chi)
 & \subset &
   r^{-n-2\beta+2} (x/r)^{2\sigma +4} e^{- 2 sr/x} C^{k+1-\lfloor n/2\rfloor + \alpha}_{x^2/r, r^{-\beta} (x/r)^{\sigma+n}e^{-sr/x}}
\\
 & = &
    C^{k+1-\lfloor n/2\rfloor + \alpha}_{x^2/r,r^{n+\beta -2} (x/r)^{-\sigma -4 + n} e^{sr/x} }
\eeal{6X16/22}
}
where $\alpha$ is any number in $(0,1)$ when $n$ is even, and $\alpha=1/2$ when $n$ is odd. (The  inclusion above is obtained by a calculation as in \cite[Equation~(B.4)]{ChDelay}  where, using the notation there, in the fifth line the  elliptic regularity estimate is replaced by the Sobolev inequality on $B(0,1/4)$ for the function $u\circ \varphi_p$.)
This guarantees differentiability of the constraints map
\blue{\referee (compare~\cite[Corollary~3.2]{bartnik:phase})}, as required for applicability of the inverse function theorem.

Smoothness of solutions is standard: away from $\partial \Omega$ this follows from elliptic regularity, while near $\partial \Omega$ this is guaranteed by the exponential decay of $h$ and its derivatives at $\partial \Omega$. The proof is now complete.
\qed

 \section{Beyond scale-invariant sets}
 \label{s19X16.1}
\newcommand{\Phihere}{\Psi}

Let $(\hyp,g)$ be an asymptotically Euclidean manifold and let  $\Omega\subset \hyp$, $r$ and $x$ be  as described at the beginning of Section~\ref{S4X16.1}.
Let $g_0$ be any smooth metric on $\hyp$ which coincides with the Euclidean metric at all large distances in the asymptotically Euclidean end. We will use the metric $g_0$ in the definitions of the functional spaces.

Suppose that $\Phihere :\hyp\to\hyp$ is a smooth diffeomorphism satisfying
\bel{19X16.1}
 \mbox{$\Psi\in   C^{k+5}_{r,r^{- 1}} $, with $\Phihere ^*g$  uniformly equivalent to $g$.}
\ee
The symbol $\chi$ will denote a cut-off function which equals the function $\chi$ defined shortly before \eq{1XI16.1} \emph{composed with $\Psi$}.

We have:

\begin{theorem}
  \label{T18X16.1}
  Under \eq{19X16.1} and the remaining hypotheses of
  Theorem~\ref{ThefulltheoremAE}, the conclusions of
  Theorem~\ref{ThefulltheoremAE}, as well as Remarks~\ref{R7X16.1}-\ref{R6X13.3}, hold  with $\Omega$ replaced by $\Phihere (\Omega)$.
\end{theorem}

\begin{remark}
 \label{R1XI16.1}
  As in our remaining analysis, the requirement that $g$ and $\hat g$ be close to each other can be achieved by translating $\Omega$ sufficiently far to the asymptotic region, or by scaling as in Section~\ref{S7X16.1}.
\end{remark}

\proof
We note that there are no KIDs on $\Phihere(\Omega)$ which tend to zero when staying away from the boundary of $\Phihere(\Omega)$. This can be proved by calculations similar to those of \cite[Proposition~2.2]{ChBeig3}, where integration over the  rays
$$
 \gamma_p=\{rp\,,\ r\in [1,\infty)\}\,,
 \
 p\in S(R_0)
 \,,
$$
is replaced by integration over $\Phihere(\gamma_p)$.

The result follows now by applying the following result to the metrics $\Phihere^* g$ and $\Phihere^* \hat g$ on $\Omega$:
\qed

\begin{theorem}
 \label{ThefulltheoremAE2}
 Under the hypotheses of Theorem~\ref{ThefulltheoremAE},
let $g_0$ be any smooth metric on $\hyp$ which coincides with the Euclidean metric in the asymptotic region. Suppose that
 $g \in C^{k+4}_{r,1}\cap C^\infty$ is uniformly equivalent to $g_0$ and has no static KIDs on $\Omega$.
For all real numbers $\sigma$ and $s>0$  and
$$
 \mbox{all smooth metrics
    $\hg $ close enough to $ g$ in $ C^{k+4}_{r,r^{- {\tilde \beta}   }}(\Omega )$}
$$
there exists on $\Omega$ a unique smooth two-covariant symmetric tensor field $ h  $ of the form 
\bel{6X16.12+}
      h
   = \phi^4\psi^2  P^*_{ g_\chi}( \delta N)\in \psi^2  \phi^2\Lpsikg{k+2}{g}
\ee
such that the metric $ g_\chi+ h   $ solves
\bel{fullcolle+}
 R
    \left[g_\chi+ h  \right]= \chi R(\hat g) +(1-\chi)
    R(g)
 \,.
\ee
The   tensor field $  h   $ vanishes at $\partial \Omega$ and can be $C^{\infty}$-extended by zero across $\partial \Omega$, leading to a smooth  metric $g_\chi+ h  $ which approaches $g$ as $r$ recedes to infinity in the asymptotic region.
\end{theorem}

\proof
Theorem~\ref{ThefulltheoremAE2}  follows from the inverse
function theorem applied to the map
\bel{28X16.1}
  r_\chi:\Lpsikg{k+4}{g} \ni \delta N
  \mapsto \psi^{-2}[(R(g_\chi +\phi^4\psi^2  P^*_{ g_\chi}( \delta N))
  - R(g_\chi)]
  \in
   \Lpsikg{k}{g}
   \,,
\ee
exactly as in the proof of Theorem~\ref{ThefulltheoremAE}.
Recall that the
 map $r_\chi$ will have $\psi^{-2}\delta R_\chi$ in its range when $g_\chi$
is close  to $g$ in $C^{k+4}_{\phi,1}$ because the inverse, $Dr_\chi(0)^{-1}$, of the linearised
map $Dr_\chi(0) $ has a bound uniform over a neighborhood of $g_\chi$ in this topology,
as follows in part from the uniformity of the constant in the weighted Poincar\'e
inequality of Section~\ref{s18IX15.6}.
\qed

\begin{remark}
  \label{R28X16.1}
An alternative proof of Theorem~\ref{T18X16.1} proceeds by
checking that Propositions~\ref{P4VIII14.1} and \ref{PL4VIII14.2} below hold with $\Omega$ replaced by $\Phihere (\Omega)$, $r$ replaced by $r_\Phihere := r\circ \Phihere^{-1} $, and $x$ replaced by $x_\Phihere :=x\circ \Phihere^{-1} $.
Since we have already seen  (see the proof of Theorem~\ref{T18X16.1}) that there are no KIDs vanishing at infinity on $\Phihere (\Omega)$,  the proof of Theorem~\ref{ThefulltheoremAE} applies verbatim.

To verify the propositions, one starts by noting that
there exists a constant $C_1$ such that we have
\bel{19X16.11}
 C^{-1}_1 r \le r_\Phihere  \le C_1 r
 \,.
\ee
This can be established by standard considerations (cf., e.g., the proof of~\cite[Equation~(30)]{ChErice}). Hence, decay rates in $\rPhi$ are equivalent to decay rates in $r$.

Let $\uPhi:= u\circ \Phihere $, and let us denote by $\mu_g$ the Riemannian measure associated with a metric $g$. We then have for functions $u$ compactly supported in $\Phihere (\Omega)$, using change-of-variables, Proposition~\ref{PL4VIII14.2}, and a trace theorem,
\bean
\lefteqn{
   \int_{\OmegaPhi  } \xPhi^{2\sigma}  \rPhi^{2\mu} e^{-s \rPhi/ \xPhi}|u| ^2 \, d\mu_g
 =    \int_{\Omega  } x^{2\sigma}  r^{2\mu} e^{-s r/x}| \uPhi|^2 \, d\mu_{\Phihere ^*g}
 }
 &&
\\
\nn
   & \le &   C \int_{\Omega  } x^{2\sigma}  r^{2\mu} e^{-s r/x}|\uPhi|^2 \, d\mu_{g}
\\
\nn
 &   \le
 &
      C^2 \bigg(
  \int_{ \{ x\ge x_0, \,  r=R \}}     |\uPhi|^2
   +
  \int_{ \{ x\ge x_0, \, r\le R \}}    |\uPhi|^2 \, d\mu_{g}
\\
\nn
  &&
 \phantom{ C\bigg(} +    \int _{\Omega   }    x^{2\sigma+4} r^{2\mu-2} e^{-sr/x}  |\nabla \uPhi|^2_g  \, d\mu_{g}
   \bigg)
\\
\nn
   & \le &   C^3  \bigg(
  \int_{ \{ x\ge x_0, \, r\le R \}}    |\uPhi|^2 \, d\mu_{\Phihere^*g}  +    \int _{\Omega   }    x^{2\sigma+4} r^{2\mu-2} e^{-sr/x}  |\nabla \uPhi|^2_{\Phihere^* g}  \, d\mu_{\Phihere^* g}
   \bigg)
\\
   & = &   C^3  \bigg(
  \int_{ \Phihere(\{ x\ge x_0, \, r\le R \})}    |u|^2 \, d\mu_{ g}  +    \int _{\Phihere(\Omega  ) }    \xPhi^{2\sigma+4} \rPhi^{2\mu-2} e^{-s\rPhi/\xPhi}  |\nabla u|^2_{g}  \, d\mu_{  g}
   \bigg)
  \,,
  \phantom{xxx}
\eeal{17X16.3}
as claimed, assuming of course  $s\ne 0$ and $ \sigma+\mu+n/2\ne 0$.

The result for tensor fields $u$ easily follows (cf., e.g., \cite[Remark~4.7]{ChDelayExotic}).

There is likewise an equivalent of \eq{9IX14.1} in the current setting.
\qedskip
\end{remark}

A trivial example of maps $\Phihere $ which satisfy the above is linear maps of $\R^n$. This, however, does not lead to a new family of sets $\Omega$ as compared to Theorem~\ref{ThefulltheoremAE}. Another example is provided by maps which are asymptotic to linear maps,  which again does not lead to any additional essential changes at large distances.

A non-trivial example is provided by \emph{logarithmic rotations}, whose action on cones or tennis-ball curves can be seen in Figure~\ref{F17X16.2}, and which we define as follows:
\begin{figure}[th]
\vspace{-.3cm}
  \centering
    \includegraphics[scale=.2]{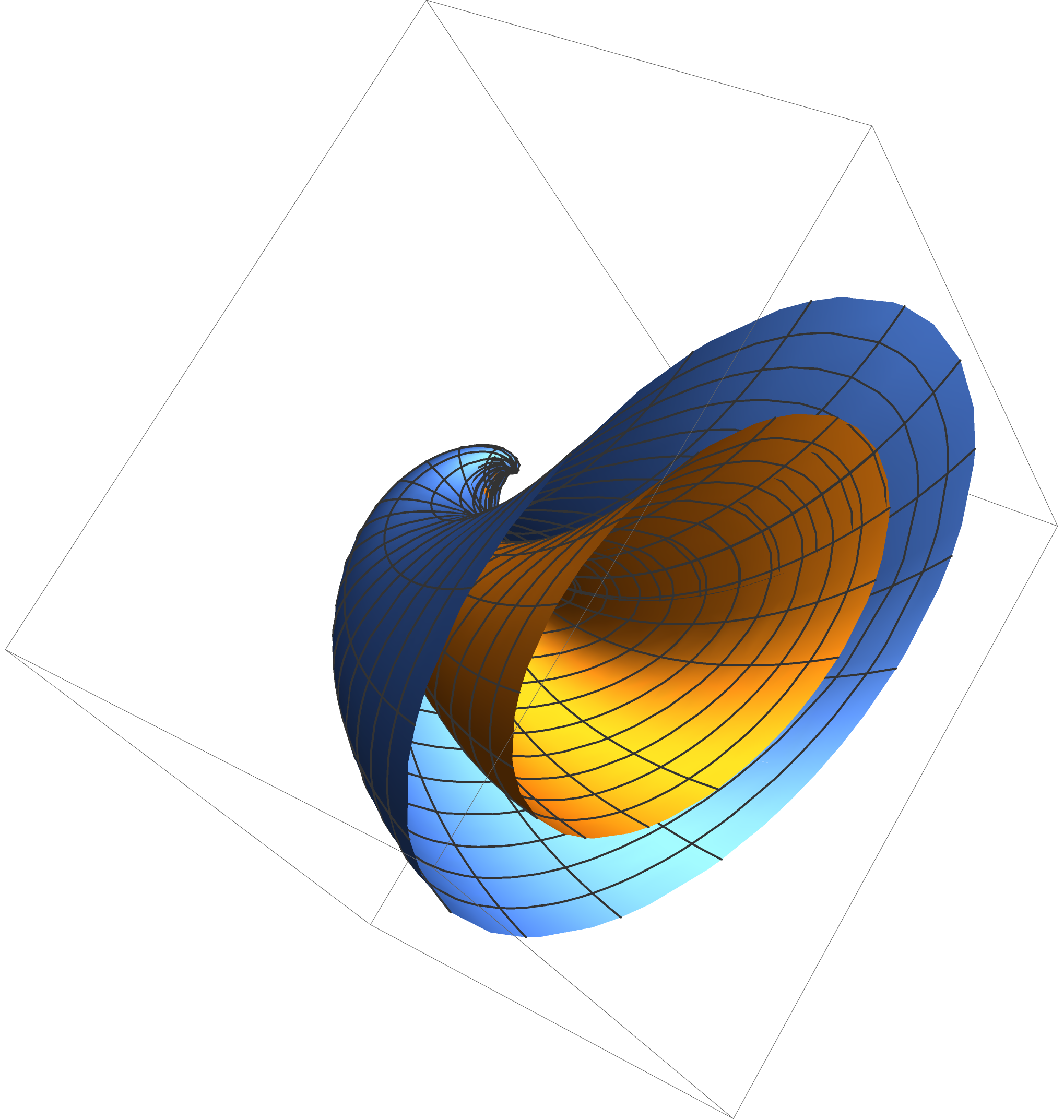}%
    \includegraphics[scale=.24]{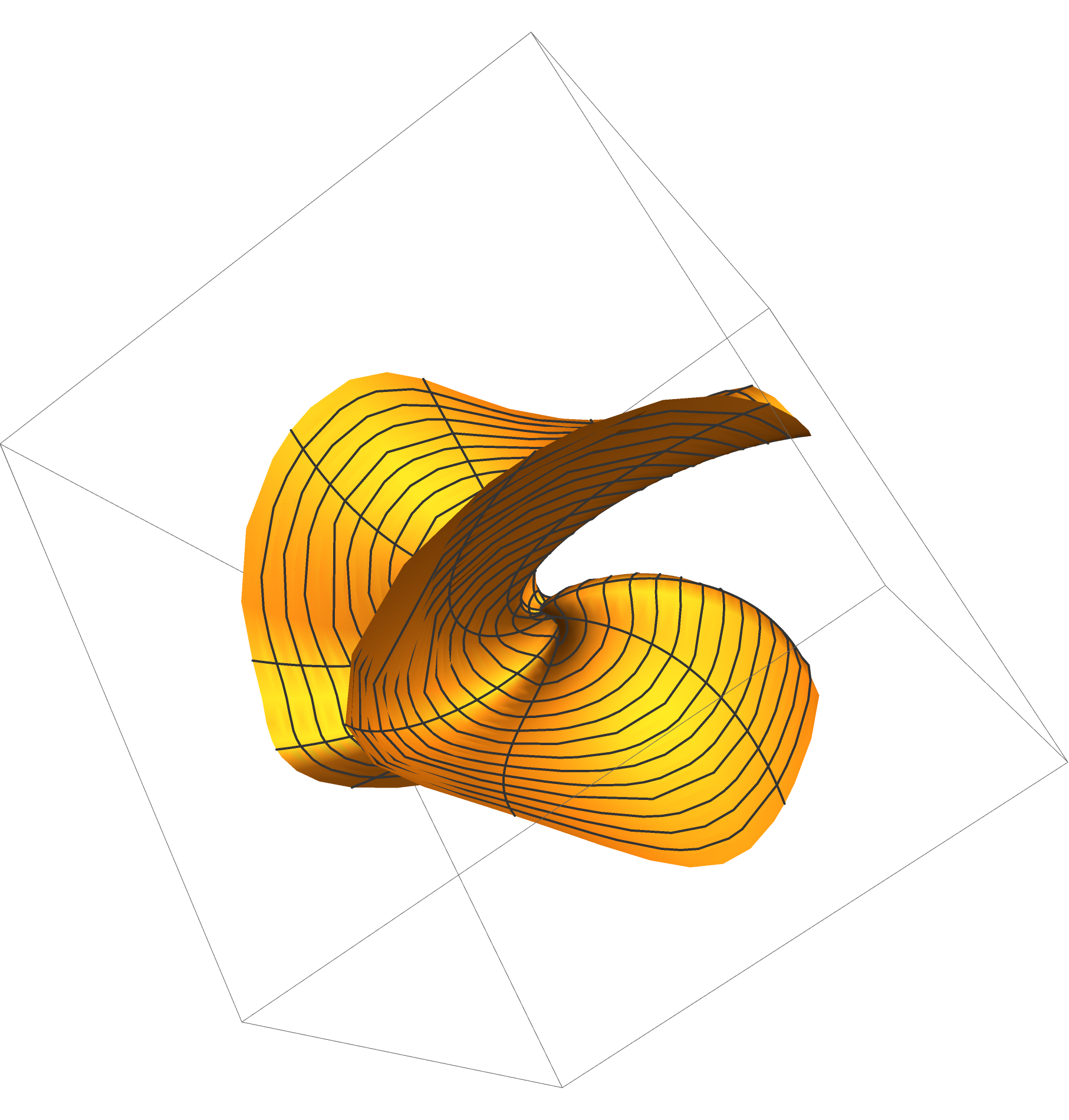}
\vspace{-.5cm}
\caption{Left figure: A  logarithmic rotation  applied to two cones; $\Omega$ is the region between the cones. Right figure: A  logarithmic rotation  applied to the ``tennis-ball curve'' (namely, the curve along which the halves of a tennis-ball are joined together). A scale-invariant thickening of the surface displayed provides the set $\Omega$.
    \label{F17X16.2}}
\end{figure}
Let $R(t)$ denote a rotation by angle $t$ around the $z$-axis. For $\rotpar\in \R$ consider the smooth map $\Phihere _\rotpar $ defined as
$$
 \R^3\setminus \ol {B(1)} \ni \vec x \mapsto R(\rotpar \ln |\vec x|) \vec x
 \in \R^3\setminus \ol {B(1)}
 \,.
$$
In spherical coordinates on $\R^3$ the map $\Phihere _\rotpar$ takes the form
$$
 (r,\theta,\varphi) \mapsto (r,\theta,\varphi+ \rotpar \ln r)
 \,.
$$
This leads to the following form of the pull-back $\Phihere _\rotpar^* \delta$ of the Euclidean metric $\delta$:
\bean
 \Phihere _\rotpar^* \delta
  & = &
   dr^2 + r^2 \Big(d\theta^2 + \sin^2\theta  \big(d \varphi + \frac {\rotpar dr} r\big)^2\Big)
\\
  & = &
   (1+\rotpar^2 \sin^2\theta) dr^2 + r^2 \Big(d\theta^2 + \sin^2\theta   d \varphi^2\Big) + 2  \rotpar r dr d\varphi
   \,.
\eeal{17X18.1}
Using
$$
 2   r dr d\varphi \le dr^2 + r^2 d\varphi^2 \le \delta
 \,
$$
we find for $|\rotpar|<1$
\bel{17X16.3+}
 (1-|\rotpar|) \delta
 \le
  \Phihere _\rotpar^* \delta \le (1+|\rotpar| +\rotpar^2) \delta \le 3 \delta
  \,,
\ee
as desired.

We note that while the image of a scale-invariant set under a logarithmic rotation   is not scale-invariant, it is invariant under a discrete scaling, which would have been enough for a direct analysis in Section~\ref{s18IX15.6} anyway.

Composing the above rotations along different axes leads to further non-trivial examples on $\R^3$.

The above generalises in an obvious way to higher dimensions by writing $\R^n=\R^3 \times \R^{n-3}$, and letting $R(t)$ be the identity on the $\R^{n-3} $ factor.

A similar calculation applies to the map
\bel{19X16.2}
 \R^3  \ni \vec x \mapsto R\Big(\rotpar \ln\sqrt{x^2+y^2+1}\Big) \vec x
 \in \R^3
 \,,
\ee
leading again to the  estimate \eq{17X16.3+} for $|\alpha|<1 $.

On $\R^{2n}$, further non-trivial maps satisfying \eq{19X16.1} can be obtained by making independent logarithmic rotations \eq{19X16.2} on each $\R^2$ factor, etc.

\section{Weighted Poincar\'e inequalities}
 \label{s18IX15.6}

We will need the following~\cite[Proposition~C.2]{ChDelay}:

\begin{Proposition}\label{prop:poinc}
Let $u$ be a $C^1$ compactly supported tensor field on a Riemannian manifold $(M,g)$, and
let $w,v$ be two $C^2$ functions defined on a neighborhood of the
support of $u$. Then
\bel{14VI14.5}
  \int_M e^{2v}|\nabla u|^2d\mu_g
    \geq  \int_M e^{2v}\left[\Delta v+\Delta w +|\nabla v|^2-|\nabla
    w|^2 \right]|u|^2d\mu_g
  \,.
\ee
\end{Proposition}

Throughout this section we assume that the metric $g$ approaches the Euclidean metric as one recedes to infinity, with no decay rate imposed. The first derivatives are required to fall-off as $o(1/r)$. Wherever second derivatives of the metric arise in a calculation, they are required to fall-off as $o(1/r^2)$, etc.
 \blue{These conditions are more than satisfied under the usual hypotheses of asymptotic flatness.}

We let $r$ be any smooth  positive function on $\Omega$ bounded away from zero which coincides with the coordinate radius on $\Omega\setminus B(R_0)$.

Our first goal is to prove Proposition~\ref{P4VIII14.1} below.   As a first step, we prove:

\begin{Lemma}
  \label{L4VIII14.1}
\blue{Let $\Omega$ be as in Section~\ref{S4X16.1}}
\referee
and let $\sigma, \mu\in \R$ be such that $\sigma\ne - 1/2$. There exist constants $c,C>0$ such that for all $C^2$ tensor fields $u$ compactly supported in
$$
 \Omega':= \Omega \cap
   \{0< x < c r \}
$$
it holds that
\bel{14VI14.5+}
  \int_\Omega x^{2\sigma+2} r^{2\mu} |\nabla u|^2d\mu_g
    \geq C\int_\Omega  x^{2\sigma } r^{2 \mu }|u|^2d\mu_g
  \,.
\ee
\end{Lemma}

\begin{remark}
  \label{R1XI16.1}
Once the inequality has been established for a metric $g $ as explained above, it is immediate that \eq{14VI14.5+} also holds, when $u$ is a function, for any metric  which is uniformly equivalent to $g $. But then it holds for tensors for all differentiable metrics uniformly equivalent to $g $ (no decay conditions on the derivatives needed)  by a standard argument;  \blue{compare~\cite[Remark~4.7]{ChDelayExotic}.}
\qed
\end{remark}

\proof
We use \eq{14VI14.5} with
\bel{14VI14.10}
 v=    (\sigma+1) \ln x + \mu \ln r
 \,.
\ee
Then
\bel{14VI14.2}
 \Delta v+ |\nabla v|^2 = \frac {\sigma+1} {  x^2}\bigg[  \sigma  |\nabla x|^2
  +  {2 } { \mu x}\frac {\nabla x \cdot \nabla r} r +   x \Delta x \bigg]
   + \frac \mu { r^2} \bigg[ (\mu-1)   {|\nabla r|^2} +   r\Delta r  \bigg]
   \,.
\ee
Letting
$$w=-\frac 12 \ln x 
$$
one similarly finds
%
\bel{14VI14.6}
 \Delta w- |\nabla w|^2 = \frac {1} {4 x^2}\bigg[    |\nabla x|^2
	- 2 x \Delta x \bigg]
   \,.
\ee
%

Given $\sigma\ne -\frac12$,
$\mu\in \R$,   
$ \epsilon>0$, and a compact set $K$, for small $x$ the dominant contribution comes from the $|\nabla x|^2/x^2$ terms. This shows that we can choose $x_0>0$ small enough so that for $0<x<x_0$ on $K$ we have
\bel{14VI14.8+}
 \Delta v+ |\nabla v|^2
 +
 \Delta w-  |\nabla w|^2
  \ge  \frac {(2\sigma+1)^2-\epsilon} {4 x^2}  |\nabla x|^2
   \,.
\ee

We consider now $p\in \Omega$ with
\bel{14VI14.7}
 \lambda := r(p) \ge 2R_0
 \,,
\ee
\blue{for an $R_0$ sufficiently large, to be determined below.}
Then the scaling transformation $\R^n\ni \vec y\mapsto 2R_0\lambda^{-1} \vec y$
 maps the set
$$ \Gamma_{\lambda}:= \Omega\cap \{ \lambda /2\le |\vec y | \le 2 \lambda \}
$$
to $\Gamma_{1}=\Omega\cap \{R_0\le |\vec y | \le 4R_0  \}$.
The associated pull-back will carry the metric $g$ on  {$\Gamma_{\lambda}$}
 to a metric $\lambda^ {-2} g_\lambda$ on {$\Gamma_{1}$}, with $g_\lambda$ approaching the flat metric on $\Gamma_{ 1}$ as $\lambda$ tends to infinity.  We can use \eq{14VI14.8+} for all rescaled metrics $g_\lambda$, with perhaps a smaller constant $x_0$ independent of $\lambda$ for $R_0$ large enough. Since  both sides are invariant   under the transformation $(x,r)\mapsto (a x, a r)$
\blue{\referee
 up to terms which are subdominant when  $R_0$ is large enough, scaling  back to  {$\Gamma_{\lambda}$}
  gives there (taking $\epsilon$ in \eqref{14VI14.8+} small enough and recalling  that  $|\nabla x| $ is close to one when $x/r$ is small)
\bel{4VIII14.1}
 \Delta v+ |\nabla v|^2
 +
 \Delta w-  |\nabla w|^2
  \ge   \frac {(\sigma+\frac 12)^2 } { 2 x^2}
   \,,
\ee
for all $0<x/r<c$ small enough, as desired.

We let $K=\overline{\Omega\cap B(0,2R_0) }$.
By compactness and \eqref{14VI14.8+} there exists a constant $c$ such that on $K$ we have for $\epsilon$ small enough
\bel{14VI14.8+a}
 \Delta v+ |\nabla v|^2
 +
 \Delta w-  |\nabla w|^2
 \ge  c\frac {(\sigma+\frac 12)^2 } {  x^2}
   \,.
\ee
This, together with \eqref{4VIII14.1}, provides the desired inequality.
}

\qedskip

We are ready to prove now:

\begin{Proposition}
  \label{P4VIII14.1}
Let $\Omega$ be as above and let $\sigma, \mu\in \R$ be such that $\sigma\ne - 1/2$ and $\sigma+\mu+n/2\ne 0$. There exist constants $x_0,C,R $ such that for all $C^2$ tensor fields $u$ compactly supported in $\Omega$
it holds that
\\
\bean
   \int_{\Omega }  x^{ 2\sigma }  r^{2\mu }  | u|^2d\mu_g
   & \le &
     C\bigg(
  \int_{ \{ x\ge x_0, \,  r=R \}}     |u|^2d\sigma_g
   +
  \int_{ \{ x\ge x_0, \, r\le R \}}    |u|^2 d\mu_g
\\
  &&
 \phantom{ C\bigg(} +    \int _{\Omega }  x^{ 2\sigma +2}  r^{2\mu } |\nabla u|^2d\mu_g
   \bigg)
  \,,
\eeal{9IX14.1}
where the symbol $d\sigma_g$  denotes the measure induced by $g$ on a submanifold of $M$.
\end{Proposition}

\proof
We will use a family of cut-off functions $\Xi_\lambda$ satisfying
\bel{8VIII14.1}
 |\phi\nabla\Xi_\lambda|\leq C/\lambda
 \,,
\ee
for $\lambda$ large. Indeed,  let
$\chi_\lambda(t)=\chi(-\ln(t)/\lambda) $
with $\chi$ a smooth non-decreasing function satisfying
$\chi(s)=0$ for $s<1$ and $\chi(s)=1$ for $s\ge 2$  \blue{(not to be confused with the function $\chi$ in \eqref {1XI16.1})}. Recall that in the current case  $\phi=x$.
 We set
 $\Xi_\lambda=\chi_\lambda(x/r)$, then
\beaa
 &
 x \nabla \Xi_\lambda
 = x \nabla \left(\chi (-\ln(x/r)/\lambda)\right)
 =
 \displaystyle
     \frac {\chi'} \lambda \left( \frac x r \nabla r -   \nabla x
      \right)
 \,,
  &
\eeaa
which is bounded on $\Omega$. Hence \eq{8VIII14.1} holds.

We note that $\Xi_\lambda=0$ if and only if $x\ge e^{-\lambda}r$, and $\Xi_\lambda =1$ if and only if $x \le e^{-2\lambda} r$.
Set
$$
 u
 =\Xi_\lambda u+(1-\Xi_\lambda)u=:v+w
 \,,
$$
thus $v$ vanishes outside  $\{0<x \le  r e^{- \lambda}\}$, and $w$ vanishes for $x \le e^{-2\lambda} r$.

Let
$$
 \psi = x^{\sigma} r^{\mu}
 \,.
$$
Applying Lemma~\ref{L4VIII14.1} to $v$ we obtain, for all $\lambda$ large enough,

\bel{11VIII14.1}
 \int_{\{0<x \le  r e^{- \lambda}\}} \psi^2\phi^2|\nabla u|^2d\mu_g
  +\frac{c_1}{\lambda^2} \int_{\{0<x \le  r e^{- \lambda}\,, \nabla \Xi_\lambda \ne 0 \}} \psi^2| u|^2d\mu_g
   \geq  c_2 \int_{\{v\ne 0\}} \psi^2| v|^2 d\mu_g
    \,.
\ee

On each half-ray we have the classical Poincar\'e inequality for tensor fields with compact support in $[R,\infty) \subset \R $: if $2\mu+2\sigma+n\ne 0$, then
\bel{8IX14.1}
 \int_R^\infty r^{2\mu + 2 \sigma+n-1} |u|^2 dr
   \le
  C_1 |u|^2(R) + C_2\int_R^\infty r^{2\mu + 2 \sigma+n+1 } |\partial_r u|^2 dr
%
     \,,
\ee
with constants $ C_2$ depending only upon $2\mu
     + 2 \sigma+n$, and $ C_1$ depending only upon $2\mu
     + 2 \sigma+n$ and $R$.
(In fact, $C_1<0$ when $2\mu+2\sigma+n>0$, but this will not be needed in our considerations.)

Integrating over the angles with $u$ replaced by $w$, this gives
\bel{14VI14.5+2}
  \int_{\{w\ne 0\}}  r^{2\mu +2\sigma +2} |\nabla w|^2 d\mu_g+
  \int_{\{w\ne 0\}\cap\{r= R\}}   | w|^2d\sigma_g
    \geq C\int_{\{w\ne 0\}\cap\{r\ge R\}}   r^{2\sigma+2 \mu }|w|^2d\mu_g
  \,.
\ee
This implies
\bean
\lefteqn{
   \int _{\{w\ne 0  \}}  r^{2\mu +2\sigma +2} |\nabla u|^2d\mu_g+ \frac{c_3}{\lambda^2}
  \int_{\{w\ne 0\,, \nabla \Xi_\lambda \ne 0 \}}  r^{2\mu +2\sigma } |u|^2d\mu_g
  }
 &&
\\
 &&
  +
  \int _{\{w\ne 0  \}\cap\{r= R\}}  |  u|^2d\sigma_g
  +
   \int _{\{w\ne 0  \}\cap\{r\le R\}}  |  u|^2d\mu_g
    \nonumber
\\
  &&   \geq
      c_4 \int _{\{w\ne 0  \}}   r^{2\sigma+2 \mu }|w|^2d\mu_g
  \,.
   \phantom{xx}
\eeal{14VI14.5+2-}

\blue{\referee
On the  support  of $w$ we have $ e^{-2\lambda}r \le x  \le C r$, and on the support of   $\nabla \Xi_\lambda$ we have  $ e^{-2\lambda}r \le x  \le e^{-\lambda}r$ } , hence
\beaa
\lefteqn{
  \frac 12 \int _{\Omega}x^{ 2\sigma }  r^{2\mu }| u|^2 d\mu_g
 \le
 \int_{\Omega} x^{ 2\sigma }  r^{2\mu }| v|^2d\mu_g+
 \int_{\Omega} x^{ 2\sigma }  r^{2\mu }| w|^2d\mu_g
 }
  &&
\\
 & \le &
 \int_{\Omega}  x^{ 2\sigma }  r^{2\mu }| v|^2d\mu_g+
 \max( \blue{C^{2\sigma}} \,, e^{{-4}\lambda \sigma}) \int_{\Omega}   r^{ 2\sigma+ 2\mu }| w|^2d\mu_g
\\
 & \le &
 \frac {C_3} 2 \bigg[ \int _{\{0<x \le  r e^{- \lambda}\}}  x^{ 2\sigma +2}  r^{2\mu } |\nabla u|^2d\mu_g
 +\frac{c_1}{\lambda^2} \int _{\{0<x \le  r e^{- \lambda}\}}  x^{ {2\sigma} }  r^{2\mu }| u|^2d\mu_g
\\
 \nonumber
&&
 +   \max( {C^{2\sigma}} \,, e^{{-4}\lambda \sigma})\bigg(
  \int _{\{w\ne 0  \}\cap\{r= R\}}  |  u|^2d\sigma_g
  +
   \int _{\{w\ne 0  \}\cap\{r\le R\}}  |  u|^2d\mu_g
\\
   &&
    +\int_{\{ x \ge  r e^{- 2\lambda}\}}  r^{2\mu +2\sigma +2} |\nabla u|^2d\mu_g+ \frac{c_3}{\lambda^2}
  \int_{\{ r e^{- 2\lambda} \le  x \le  r e^{- \lambda}\}}   r^{2\mu +2\sigma } |u|^2d\mu_g
  \bigg)\bigg]
  \,.
\eeaa
\blue{Without loss of generality we can assume that $c_3\ge 1$.}
For $\lambda \ge 1$ we can rewrite this as
\beaa
\lefteqn{
   \int_{\Omega} x^{ 2\sigma }  r^{2\mu }\left(1 - C_3 \frac{c_1}{\lambda^2}\right) | u|^2d\mu_g
  }
  &&
\\
 & \le &
    C_3 \bigg[ c_3
     \max( \blue{C^{2\sigma}} \,, e^{{-4}\lambda \sigma})
     \bigg(
  \int_{\{ r e^{- 2\lambda} \le  x \le  r e^{- \lambda}\}}    r^{2\mu +2\sigma } |u|^2d\mu_g
  +
  \int _{\{w\ne 0  \}\cap\{r= R\}}  |  u|^2d\sigma_g
\\
&&
  +
   \int _{\{w\ne 0  \}\cap\{r\le R\}}  |  u|^2d\mu_g
 +    \int_{\{ x \ge  r e^{- 2\lambda}\}}  r^{2\mu +2\sigma +2} |\nabla u|^2d\mu_g
 \bigg)
\\
  &&
  +
    \int _{\{0<x \le  r e^{- \lambda}\}}  x^{ 2\sigma +2}  r^{2\mu } |\nabla u|^2d\mu_g
  \bigg]
  \,.
\eeaa
Choosing $\lambda$ large enough we obtain, keeping in mind that {$x/r$ is uniformly bounded  on $\Omega$,  while $r \le  xe^{2\lambda} $ in the \blue{second-to-last}
\referee
 integral so that $r^{2\sigma+2}$ is equivalent to $  x  ^{2\sigma+2}$ there,}
\bean
\lefteqn{
   \int _{\Omega} x^{ 2\sigma }  r^{2\mu }  | u|^2d\mu_g
  \le
      C_4 (\lambda)
      \bigg(
  \int_{\{  r e^{- 2\lambda} \le  x \le  r e^{- \lambda} \}}   r^{2\mu +2\sigma } |u|^2d\mu_g
  +
  \int _{\{w\ne 0  \}\cap\{r= R\}}  |  u|^2d\sigma_g
}
 &&
\\
 &&
 \phantom{xxxxxxxxxx}
  +
   \int _{\{w\ne 0  \}\cap\{r\le R\}}  |  u|^2
  +     \int _M  x^{ 2\sigma +2}  r^{2\mu } |\nabla u|^2d\mu_g
  \bigg)
  \,.
\eeal{8IX14.4}
Integrating \eq{8IX14.1} over the angles gives
\bean
 \lefteqn{
  \int_{\{  r e^{- 2\lambda} \le  x \le  r e^{- \lambda}\,, r\ge R \}}   r^{2\mu +2\sigma } |u|^2d\mu_g
  }
  &&
\\
 \nonumber
 &&
   \le
 C(R)
  \int_{\{  r e^{- 2\lambda} \le  x \le  r e^{- \lambda}\,, r\ge R  \}}   r^{2\mu +2\sigma +2} |\nabla u|^2d\mu_g
   +
 C_1(R)
  \int_{\{  r e^{- 2\lambda} \le  x \le  r e^{- \lambda}\,, r= R  \}} | u|^2d\sigma_g
\\
 \nonumber
 &&
   \le
 C'(R,\lambda)
  \int_{\{  r e^{- 2\lambda} \le  x \le  r e^{- \lambda} \}}  x^{2\sigma +2} r^{2\mu } |\nabla u|^2d\mu_g
   +
 C_1(R)
  \int_{\{  r e^{- 2\lambda} \le  x \le  r e^{- \lambda}\,, r= R  \}} | u|^2d\sigma_g
\\
 &&
   \le
 C'(R,\lambda)
  \int_{M}  x^{2\sigma +2} r^{2\mu } |\nabla u|^2d\mu_g
   +
 C_1(R)
  \int_{\{  r e^{- 2\lambda} \le  x \le  r e^{- \lambda}\,, r= R  \}} | u|^2d\sigma_g
 \,.
\eeal{8IX14.2}
Inserting this into the first line of \eq{8IX14.4} gives \eq{9IX14.1} \,.
\qedskip

We note a standard consequence of Proposition~\ref{P4VIII14.1}:

\begin{Corollary}
  \label{C4VIII14.1}
Let $\Omega$ be as above and let $\sigma, \mu\in \R$ be such that $\sigma\ne - 1/2$ and $\sigma+\mu+n/2\ne 0$. Let $\mathring H$ denote the completion of $C^2$    compactly supported tensor fields in $\Omega$ with respect to the norm
$$
 \|u\|_{\mathring H} ^2=   \int_{\Omega }  x^{ 2\sigma }  r^{2\mu }  | u|^2d\mu_g
 +\int _{\Omega }  x^{ 2\sigma +2}  r^{2\mu } |\nabla u|^2d\mu_g
    \,.
$$
Suppose that  $ \mathring H$ contains  a closed subspace $X$ transversal to the space $\{u\in \mathring H: \nabla u = 0\}$. Then there exists a constant $C$ such that for all $u\in X$ we have
\bea
 \int_{\Omega }  x^{ 2\sigma }  r^{2\mu }  | u|^2d\mu_g
   & \le &
     C  \int _{\Omega }  x^{ 2\sigma +2}  r^{2\mu } |\nabla u|^2d\mu_g
  \,.
\eeal{11IX14.130}
\end{Corollary}

\proof
As already mentioned, the result is standard, we give the proof for completeness.

Suppose that  \eq{11IX14.130} is wrong, then there exists a sequence of $C^2$ compactly supported tensor fields $u_n$ such that
\bea
 \int_{\Omega }  x^{ 2\sigma }  r^{2\mu }  | u_n|^2d\mu_g
   & \ge &
     n  \int _{\Omega }  x^{ 2\sigma +2}  r^{2\mu } |\nabla u_n|^2d\mu_g
  \,.
\eeal{11IX14.2}
We can normalize the sequence $u_n$ so that
\bea
  \int_{ \{ x\ge x_0, \,  r=R \}}     |u_n|^2d\sigma_g
   +
  \int_{ \{ x\ge x_0, \, r\le R \}}    |u_n|^2d\mu_g
  =1
  \,.
\eeal{9IX14.1+c}
\Eq{9IX14.1} implies
\bean
   \int_{\Omega }  x^{ 2\sigma }  r^{2\mu }  | u_n|^2d\mu_g
   & \le &
     C\bigg(1 +
  \frac 1 n
  \int_{\Omega }  x^{ 2\sigma }  r^{2\mu }  | u_n|^2d\mu_g
   \bigg)
  \,.
\eeal{9IX14.1+a}
Thus, for $n\ge 2C$,
\bea
   \int_{\Omega }  x^{ 2\sigma }  r^{2\mu }  | u_n|^2d\mu_g
   & \le &
     2 C
  \,.
\eeal{9IX14.1+b}
\Eq{11IX14.2} gives now
\bea
 \int _{\{x\ge x_0,\, r\le R  \}}   |\nabla u_n|^2d\mu_g
 \le  C_1
 \int _{\Omega }  x^{ 2\sigma +2}  r^{2\mu } |\nabla u_n|^2d\mu_g
 \le \frac{2 C C_1}n
  \,.
\eeal{11IX14.1+ }
Let $K= \Omega\cap \{x> x_0,\, r< R \}$.
Compactness of the embeddings
$$
W^{1,2}(K)
 \subset
 L^2(K)\,,
 \quad
 W^{1,2}(K)
 \subset
 L^2(\{r=R\}\cap K)\,,
$$
implies that $\{u_n\}$ contains a subsequence, still denoted by  $\{u_n\}$, which is Cauchy both in $L^2(K)$ and in $L^2(\{r=R\}\cap K)$. \Eq{9IX14.1} applied to $u_n-u_m$ shows that $\{u_n\}$ is Cauchy in $\mathring H$. The limit is a non-trivial tensor field satisfying $\nabla u =0$, which contradicts the fact that zero is the only such tensor field in $X$.
\qedskip

We have the exponentially-weighted version of the above:

\begin{Proposition}
  \label{PL4VIII14.2}
Let $\Omega$ be as above and let $s, \sigma,  \mu\in \R$ satisfy $s\ne 0$ and $\beta\equiv {-\sigma-\mu-n/2}\ne 0$.
 There exist constants $x_0,C,R $ such that for all $C^2$ tensor fields $u$ compactly supported in $\Omega$
it holds that
\bean
   \int_{\Omega } x^{2\sigma} r^{2\mu} e^{-sr/x}|u|^2d\mu_g
   & \le &
     C\bigg(
  \int_{ \{ x\ge x_0, \,  r=R \}}     |u|^2d\sigma_g
   +
  \int_{ \{ x\ge x_0, \, r\le R \}}    |u|^2d\mu_g
\\
  &&
 \phantom{ C\bigg(} +    \int _{\Omega }    x^{2\sigma+4} r^{2\mu-2} e^{-sr/x}  |\nabla u|^2d\mu_g
   \bigg)
  \,.
\eeal{9IX14.110}
\end{Proposition}

\proof
The proof is a direct repetition of that of Proposition~\ref{P4VIII14.1}, with $\psi= x^\sigma r^\mu e^{-sr/(2x)}$, $\phi =x^2/r$, and where instead of Lemma~\ref{L4VIII14.1} the following Lemma is used:
\qedskip

\begin{Lemma}
  \label{L4VIII14.2}
Let $\Omega$ be as above and let $s, \mu\in \R$ be such that $s\ne 0$.
  There exist constants $c,C>0$ such that for all $C^2$ tensor fields $u$ compactly supported in
$$
\Omega':= \Omega \cap
   \{0<  x\ < c r \}
$$
it holds that
\bel{14VI14.5+2*}
  \int_M x^{4+2\sigma} r^{2\mu-2} e^{-sr/x} |\nabla u|^2d\mu_g
    \geq C\int_M x^{2\sigma }r^{2\mu} e^{-sr/x}|u|^2d\mu_g
  \,.
\ee
\end{Lemma}

\proof
We use  again \eq{14VI14.5}, with $w=0$ and
\bel{14VI14.10x}
 v=  -\frac {sr} x + (2 +\sigma) \ln x + (\mu-1) \ln r
 \,.
\ee
On any compact set one finds, for $0<x/r$ small enough,
\bel{14VI14.2x}
 \Delta v+ |\nabla v|^2 =
 \big(
  \frac {s^2r^2}{x^4} + O(r x^{-3})
  + O( x^{-2})
   + O(r^{-2})
   \big)
    |\nabla x |^2
   \,.
\ee
\blue{\referee
A scaling argument as in the proof of Lemma~\ref{L4VIII14.1} shows that
} for $s\ne 0$ and for $0<x<cr$  with $c$ small enough we have
\bel{14VI14.2x3}
 \Delta v+ |\nabla v|^2 \ge \hat C  \frac {r^2}{x^4}
   \,,
\ee
for some constant $\hat C >0$. This,
together with \eq{14VI14.5} leads to \eq{14VI14.5+2*}.
\qedskip

\section{Applications}
 \label{S7X16.1}

In this section we wish to show how to arrange things so that all conditions of Theorem~\ref{ThefulltheoremAE} are met.

Consider, then, two smooth metrics $g, \hat g \in M_{\delta+C^{k+4}_{r,r^{\epsilon}}}\cap C^\infty$, $i=1,2$, for some $\epsilon>0$ and $k>n/2$. Let $\Omega$, $\Omega_S$ be as in Section~\ref{S4X16.1}, cf. \eq{7X16.24}-\eq{7X16.25}.

Let $\psi_{\vec y}:\R^n\to\R^n$ denote the   translation by $\vec y$,
$$
 \phi_{\vec y}(\vec x):= \vec x + \vec y
 \,.
$$
When $|\vec y|$ tends to infinity the three tensor fields $\phi_{\vec y}^* g -\delta$, $\phi_{\vec y}^* \hat g -\delta$, and $\phi_{\vec y}^*(g-\hat g)$   approach   zero in $ C^{k+4}_{r,r^{-\epsilon  }}(\Omega )$. Thus Theorem~\ref{ThefulltheoremAE} applies for all $|\vec y|$ large enough, and provides a gluing of $\phi_{\vec y}^* g $ with $\phi_{\vec y}^* \hat g$. Taking $\Omega$ to be a translation of the asymptotic cone as in \eq{7X16.31} proves Theorem~\ref{T11IX16.1}.

An alternative construction proceeds as follows, assuming that $\Omega\cap B(0,1)=\emptyset$: For $\lambda \ge 1$ let $\psi_\lambda: \R^n\to \R^n$ denote the scaling map, $\psi_\lambda(\vec x) := \lambda \vec x$. Set
\bel{7X16.26}
 g_\lambda = \lambda^{-2} \psi_\lambda^* g
 \,,
  \qquad
 \hat g_\lambda = \lambda^{-2} \psi_\lambda^* \hat g
 \,.
\ee
Then both $g_\lambda$ and $\hat g_\lambda$ tend to $\delta$ in $ C^{k+4}_{r,r^{-\epsilon }}(\Omega )$ as $\lambda$ tends to infinity.  Hence Theorem~\ref{ThefulltheoremAE} applies for all $\lambda$ large enough. Scaling back the glued  metric provides a gluing of $g$ and $\hat g$ along a rescaled set $\Omega$.

Consider, finally, a collection $\Omega_i$, $i=1,\ldots,N$,  of disjoint sets satisfying the requirements set forth at the beginning of Section~\ref{S4X16.1}, possibly after translations. The construction just given can be repeated $N$-times to obtain a gluing of $g$ and $\hat g$ across $\cup \Omega_i$.

\bigskip

\noindent{\sc Acknowledgements:} Supported in part by the
Austrian Research Fund (FWF), Project P 29517-N16, and by the grant ANR-17-CE40-0034 of the French National Research Agency ANR (project CCEM).

\appendix

\bibliographystyle{amsplain} \providecommand{\bysame}{\leavevmode\hbox to3em{\hrulefill}\thinspace}
\providecommand{\MR}{\relax\ifhmode\unskip\space\fi MR }
\providecommand{\MRhref}[2]{%
  \href{http://www.ams.org/mathscinet-getitem?mr=#1}{#2}
}
\providecommand{\href}[2]{#2}

\end{document}